\documentclass[12pt,a4paper]{amsart}

\usepackage{fullpage}
\usepackage[pdftex]{epsfig}

\usepackage{amssymb}  
\usepackage{latexsym} 
\usepackage{comment}
\usepackage{color}
\usepackage{enumerate}
\usepackage[all]{xy}

\usepackage{charter,euler} 
\DeclareSymbolFont{bchoperators}{T1}{bch}{m}{n}
\SetSymbolFont{bchoperators}{bold}{T1}{bch}{b}{n}
\makeatletter
\renewcommand{\operator@font}{\mathgroup\symbchoperators}
\makeatother

\DeclareFontEncoding{OT2}{}{} 
\newcommand{\textcyr}[1]{%
 {\fontencoding{OT2}\fontfamily{wncyr}\fontseries{m}\fontshape{n}\selectfont #1}}
\newcommand{\Sha}{{\mbox{\textcyr{Sh}}}}

\usepackage{titlesec} 
\titleformat{\section}{\normalfont\bfseries\filcenter}{\thesection}{1em}{}
\titleformat{\subsection}{\normalfont\bfseries}{\thesubsection}{1em}{}
\titleformat{\subsubsection}{\normalfont\bfseries}{\thesubsubsection}{1em}{}

\newcommand{\Z}{{\mathbb Z}}
\newcommand{\Q}{{\mathbb Q}}
\newcommand{\R}{{\mathbb R}}
\newcommand{\F}{{\mathbb F}}

\newcommand{\PP}{{\mathbb P}}

\newcommand{\CK}{{\mathcal K}}
\newcommand{\CO}{{\mathcal O}}

\newcommand{\fp}{{\mathfrak p}}
\newcommand{\To}{\longrightarrow}

\newcommand{\Sel}{\operatorname{Sel}}

\newcommand{\tors}{{\operatorname{tors}}}
\newcommand{\Pic}{\operatorname{Pic}}
\newcommand{\im}{\operatorname{im}}

\newenvironment{steplist}[1]%
   {\begin{list}{}{\settowidth{\labelwidth}{#1}%
                   \setlength{\leftmargin}{\labelwidth}%
                   \addtolength{\leftmargin}{\labelsep}}}%
   {\end{list}}

\newtheorem{Theorem}{Theorem}[section]
\newtheorem{Lemma}[Theorem]{Lemma}
\newtheorem{Proposition}[Theorem]{Proposition}

\newtheorem{Conjecture}[Theorem]{Conjecture}

\theoremstyle{definition}

\newtheorem{Examples}[Theorem]{Examples}
\newtheorem{Remark}[Theorem]{Remark}

\newtheorem{Algorithm}[Theorem]{Algorithm}

\numberwithin{equation}{section}

\definecolor{darkgreen}{rgb}{0,0.5,0}
\definecolor{rem}{rgb}{0.8,0,0}
\definecolor{new}{rgb}{0.3,0.1,0.9}
\definecolor{reply}{rgb}{0,0,0.8}
\definecolor{gray}{gray}{0.7}
\usepackage[
        colorlinks, citecolor=darkgreen,
        backref,
        pdfauthor={Michael Stoll},
        pdftitle={Determining the rational points on a curve of genus 2 and Mordell-Weil rank 1},
]{hyperref}
\usepackage[alphabetic,backrefs,lite]{amsrefs} 

\begin{document}

\title{Determining the rational points on a curve of genus 2 \\ and Mordell-Weil rank 1}

\author{Michael Stoll}
\address{Mathematisches Institut,
         Universität Bayreuth,
         95440 Bayreuth, Germany}
\email{Michael.Stoll@uni-bayreuth.de}
\urladdr{\url{http://www.mathe2.uni-bayreuth.de/stoll/}}

\date{\today}


\begin{abstract}
  We explain how one can efficiently determine the (finite) set of rational points
  on a curve of genus~$2$ over~$\Q$ with Jacobian variety~$J$, given a point
  $P \in J(\Q)$ generating a subgroup of finite index in~$J(\Q)$.
\end{abstract}

\maketitle


\section{Introduction} 

Let $C$ be a (``nice'': smooth, projective and geometrically irreducible)
algebraic curve
of genus~$2$ over~$\Q$. Concretely, $C$ can be specified by a polynomial
$f \in \Z[x]$, squarefree and of degree $5$ or~$6$; then $C$ is the smooth
projective model of the affine plane curve given by
\[ y^2 = f(x) \,. \]
We denote the Jacobian variety of~$C$ by~$J$.

By Faltings's Theorem \cite{Faltings,FaltingsErratum}, whose statement was first
conjectured by Mordell \cite{Mordell}, we know that
the set~$C(\Q)$ of $\Q$-rational points on any nice algebraic curve~$C$ over~$\Q$
of genus at lest~$2$ is finite. (More generally, the corresponding result
holds for curves over any algebraic number field.) However, all known
proofs (see~\cite{BombieriGubler}*{Chapter~11} and~\cite{LawrenceVenkatesh}
for two alternative proofs)
of this theorem are \emph{infeffective:} they do not give us
an algorithm (not even a terribly inefficient one) that would determine
this finite set in any concrete case. (It should be noted that there
is a general algorithm that would compute the set of rational points on
a curve of higher genus, whose termination is conditional on the
Hodge, Tate and Fontaine-Mazur conjectures; see~\cite{AlpogeLawrence}.
However, this algorithm is so complex that it is unlikely to finish
for any concrete example before the end of the universe.)

So it is an interesting problem to determine this set~$C(\Q)$ explicitly and
in reasonable time. The setting of curves of genus~$2$ over~$\Q$ is
the simplest possible case; it is therefore natural to consider it.

The purpose of this note is to explain how $C(\Q)$ can be computed if we know a
point $P \in J(\Q)$ such that $P$ generates a subgroup of finite index in~$J(\Q)$.
It is known (due to Weil~\cite{Weil}) that $J(\Q)$ is a finitely generated abelian group;
the assumption above then implies that its rank is at most~$1$.
In particular, the rank of~$J(\Q)$ is strictly
less than the genus of the curve. In this setting, Chabauty~\cite{Chabauty}
gave a proof of Mordell's Conjecture some forty years before Faltings
proved it in general. Coleman~\cite{Coleman} showed how Chabauty's approach
can be used to obtain quite good explicit \emph{bounds} on the \emph{number} of
rational points. This approach can also be used to \emph{find} the set of rational
points. This was first done in~\cite{FPS}; further expositions can be found
in~\cite{McCP,Stoll2006a}, and another concrete application, e.g., in~\cite{Stoll2008a}.

We note that a rational
point $P \in J(\Q) \setminus \{0\}$ can be specified by its \emph{Mumford representation};
this is a pair of polynomials $a, b \in \Q[x]$, with $a$ monic of degree
at most~$2$, such that $P$ is in the linear equivalence class of $D - D_\infty$,
where $D$ is cut out by $a(x) = 0$, $y = b(x)$ (when $\deg(a) < 2$, one has
to modify this slightly to include points at infinity) and $D_\infty$ is
the polar divisor of~$x$.

Our main result in this note is the following.

\begin{Theorem} \label{thm:main}
  There is an algorithm that
  \begin{enumerate}[\upshape (1)]
    \item takes as input a polynomial $f \in \Z[x]$ as above and (the Mumford
          representation of) a point $P \in J(\Q)$;
    \item terminates assuming that Conjectures~\ref{conj:MWS} and~\ref{conj:diff} hold for~$C$;
    \item upon termination returns a subset of~$C(\Q)$ that is all of~$C(\Q)$
          when $(J(\Q) : \langle P \rangle) < \infty$.
  \end{enumerate}
\end{Theorem}

This algorithm has been implemented in the Magma computer algebra system~\cite{Magma};
it is included in the package files shipped with Magma as of version~2.29.
In practice, it terminates within a few seconds on most inputs of reasonable size.

Here is a rough outline of the algorithm. We assume that $(J(\Q) : \langle P \rangle) < \infty$.
\begin{enumerate}[1.]
  \item \label{alg:els}
        Check whether $C$ has points everywhere locally (i.e., over all completions of~$\Q$).
        If not, then $C(\Q) = \emptyset$; Stop.
  \item \label{alg:tors}
        Determine the torsion subgroup~$J(\Q)_\tors$ of~$J(\Q)$.
  \item \label{alg:rk0}
        If $P$ is of finite order, then $J(\Q)$ is torsion. We can determine~$C(\Q)$
        via its image in~$J(\Q)$ under $Q \mapsto [Q - \iota(Q)]$
        ($\iota \colon (x,y) \mapsto (x,-y)$ is the hyperelliptic involution); Stop.
  \item \label{alg:pic1}
        Check whether $C$ has rational divisors of odd degree.
        If not, then $C(\Q) = \emptyset$; Stop.
  \item \label{alg:mws}
        Use a rational divisor of odd degree to embed $C$ into~$J$ and run a
        ``Mordell-Weil sieve + Chabauty'' computation.
        This will terminate in practice (and in theory subject to Conjectures
        \ref{conj:MWS} and~\ref{conj:diff}) and in this case return~$C(\Q)$.
\end{enumerate}

Before we go into details regarding the steps of the above algorithms,
we give  in Section~\ref{S:MW} a short summary of how we can obtain generators
of a finite-index subgroup of the Mordell-Weil group~$J(\Q)$.
This is essential in order
to obtain the necessary input, namely, to show that the rank~$r$ is at most~$1$ and
that we have a point~$P$ that generates a subgroup of rank~$r$.

The computation in Step~\ref{alg:els} is essentially standard.
Section~\ref{S:ELS} discusses an efficient implementation for hyperelliptic curves.

An algorithm for Step~\ref{alg:tors} is described in~\cite{Stoll1999}*{\S\,11}
and is available in Magma as \texttt{TorsionSubgroup}($J$).

Step~\ref{alg:rk0} has been available in Magma for a while as \texttt{Chabauty0}($J$).

Step~\ref{alg:pic1} is new. See Section~\ref{S:Pic1} below for a description.

The basic algorithm for Step~\ref{alg:mws} is described in~\cite{BruinStoll2010}.
What is new here is a complete treatment of the case when the Jacobian splits.
See Section~\ref{S:MWS} for details.


\section{Computing the Mordell-Weil group} \label{S:MW}

Let $C$ be a nice curve of genus $g \ge 2$ over a number field~$K$, with Jacobian variety~$J$.
By the (Mordell-)Weil Theorem~\cite{Weil}, the group~$J(K)$ of $K$-rational points on~$J$
is a finitely generated abelian group; it therefore has a finite torsion subgroup~$J(K)_\tors$,
and the quotient $J(K)/J(K)_\tors$ is a free abelian group of some rank $r \in \Z_{\ge 0}$,
which is the \emph{rank} of~$J(K)$.

It is then an interesting problem and also an important step in the solution
of various other problems related to the arithmetic of~$C$ to determine~$J(K)_\tors$,
the rank~$r$, and explicit points $P_1, \ldots, P_r \in J(K)$ whose images are
free generators of the quotient~$J(K)/J(K)_\tors$.

The usual approach is as follows.
\begin{enumerate}[1.]
  \item Search for points in~$J(K)$ and check for relations between them.
        This will give a lower bound on~$r$.
  \item Compute a Selmer group (and use perhaps further methods) to obtain
        an upper bound on~$r$.
  \item \label{MW:step3}
        If both bounds agree, $r$ is determined, and we know points that
        generate a subgroup of finite index. If the bounds do not agree,
        go back and try to improve one or both of the bounds.
  \item Saturate the known finite-index subgroup to obtain generators
        of the full group.
\end{enumerate}
For any $n \ge 2$, there is the \emph{$n$-Selmer group}~$\Sel_n(J/K)$ of~$J$ over~$K$.
It is a finite $\Z/n\Z$-module sitting in an exact sequence
\[ 0 \To J(K)/nJ(K) \To \Sel_n(J/K) \To \Sha(J/K)[n] \To 0 \,, \]
where $\Sha(J/K)$ is the Tate-Shafarevich group of~$J$ over~$K$.
When $n = p$ is a prime, then $\Sel_p(J/K)$ is a finite-dimensional $\F_p$-vector
space, and we have
\[ \dim_{\F_p} \Sel_p(J/K) = r + \dim_{\F_p} J(K)[p] + \dim_{\F_p} \Sha(J/K)[p] \,, \]
so (assuming we can determine $J(K)[p]$) we obtain a bound on~$r$ from it.
At least in principle, Selmer groups are computable. In practice, we can
usually compute the $2$-Selmer group of the Jacobian of a \emph{hyperelliptic}
curve $C \colon y^2 = f(x)$, when the number fields generated by the roots of~$f$
over~$K$ are not ``too large''. (The main bottleneck is the computation of the
class groups of these number fields. In practice, one frequently works assuming
GRH to speed up this part of the computation.) See~\cite{Cassels,CasselsFlynn}
for the case of genus~$2$ and~\cite{PoonenSchaefer,Schaefer,BPS} for a more
general setting. The algorithm for the computation of $2$-Selmer groups of
hyperelliptic Jacobians is described in detail in~\cite{Stoll2001}
(and implemented in Magma).

There are two (not mutually exclusive) ways in which Step~\ref{MW:step3} above can fail:
there may be nontrivial $p$-torsion in~$\Sha(J/K)$; then the upper bound
will not be tight, or we missed some generators in our search, because they
are too large, so the lower bound is not tight. We now assume that $p = 2$.
We can determine what the parity of the $\F_2$-dimension of~$\Sha(J/K)[2]$
should be~\cite{PoonenStoll1999b}; this will be correct assuming that the
$2$-primary part of~$\Sha(J/K)$ is finite. Note that it is generally conjectured
that~$\Sha(J/K)$ is finite. If the parity is odd, then we can subtract~$1$
from the upper bound. Then if the upper bound is not tight, the difference
with the lower bound must be at least~$2$, so this can give an indication
whether it is worth while trying to improve the upper bound. There are various
ways in which this can be attempted. One is to consider other abelian varieties~$A$
that are isogenous over~$K$ (by an isogeny of $2$-power degree) to~$J$; then
$A$ and~$J$ have the same rank, but it is possible that the upper bound
for the rank of~$A$ that we obtain from its $2$-Selmer group is smaller
than that obtained for~$J$. The Magma implementation (for genus~$2$ Jacobians
over~$\Q$) looks at all $2$-power isogenous Jacobians, products of elliptic
curves and Weil restrictions of elliptic curves over quadratic number fields
for this purpose. Another possibility is to try to \emph{visualize} nontrivial
elements of~$\Sha(J/K)[2]$ in another abelian variety that shares some part
of the $2$-torsion Galois module with~$J$; see~\cite{CremonaMazur,BruinFlynn}.
In our implementation we work with quadratic twists of~$J$. Finally, one can
try to compute the Cassels-Tate pairing~\cite{Cassels1962,Tate} on the $2$-Selmer group,
for which there is a recent algorithm by Fisher and Yan~\cite{FisherYan}
(this algorithm is not yet part of the Magma distribution and so is not
currently used in our implementation).

To search for rational points on the (genus~$2$) Jacobian (over~$\Q$),
one can use the \texttt{j-points} program written by the author of this note,
which is included in Magma. It uses a quadratic sieve to find rational
points on the associated Kummer surface that lift to the Jacobian.
This approach has cubic complexity in the bound for the multiplicative
height of the points we want to find; it is therefore not of much use
when we want to find points of larger height. Instead, we can construct
the $2$-covering spaces associated to the elements of the Selmer group
that are not yet hit by the points we found and search for rational
points on them. The advantage of this approach is that the points have
smaller height on the covering spaces, but these spaces are geometrically
more complicated than the Kummer surface.

For the saturation step we first find a bound on the index of the subgroup~$G$
we have found in~$J(\Q)$. We then check for each prime~$p$ below this bound
whether $G$ is already $p$-saturated, i.e., whether the natural map
$G/pG \to J(\Q)/pJ(\Q)$ is an isomorphism. This can be achieved by considering
the image of~$G$ in~$J(\F_\ell)$ for various primes~$\ell$ such that
$p$ divides~$\#J(\F_\ell)$. Even when $G$ is not $p$-saturated, this will
give strong restrictions on the elements of~$G$ that might be divisible by~$p$
in~$J(\Q)$; we can then check enough of these elements individually.
Note that the input to the algorithm in Theorem~\ref{thm:main} need not be saturated. However,
in the course of running the algorithm, we will need to saturate the known
subgroup at certain primes, and so we need a procedure that does the $p$-saturation
for a given prime~$p$.

A function that uses the approach described here to try to determine generators
of the Mordell-Weil group~$J(\Q)$ when $J$ is the Jacobian variety of a
curve of genus~$2$ over~$\Q$ is available as \texttt{MordellWeilGroupGenus2}
(or also as part of the more general \texttt{MordellWeilGroup} function) in~Magma.


\section{Testing for local points} \label{S:ELS}

In this section we assume that $C$ is a hyperelliptic curve of genus~$g$ over
a number field~$K$ with affine model of the form $y^2 = f(x)$,
where $f$ is a squarefree polynomial with coefficients in~$K$ and
$\deg f \in \{2g+1, 2g+2\}$. We assume that $g \ge 2$ in the following.

We do not claim that this section contains any original results; the material
is certainly well-known to the experts. We use the shorthand ``ELS''
for ``everywhere locally soluble'', i.e., the statement that $C(K_v) \neq \emptyset$
for all places~$v$ of~$K$.

The algorithm has two parts:
\begin{enumerate}[(1)]
  \item Determine a set~$S$ of places of~$K$ such that $C(K_v) \neq \emptyset$
        for all $v \notin S$ can be shown without computation at individual places,
        with $S$ as small as possible.
  \item For each $v \in S$, decide whether $C(K_v) = \emptyset$ or not.
\end{enumerate}

We note that when $\deg(f) = 2g+1$ is odd, then there is a $K$-rational
point at infinity, and so in particular, the curve is ELS. So we will
assume that $\deg(f) = 2g+2$.

We begin with the first part. Recall the following two well-known facts.

\begin{Lemma}[Lower Weil bound] \label{L:Weil}
  If $C$ is a nice curve of genus~$g$ over a finite field~$\F_q$, then
  \[ \#C(\F_q) \ge q - 2g \sqrt{q} + 1 \,. \]
\end{Lemma}

\begin{Lemma}[Hensel lifting] \label{L:Hensel}
  If $C$ is a nice curve over a $p$-adic field~$K$ with residue class field~$k$
  and there is a smooth $k$-rational point~$P$ on the reduction of~$C$, then
  $C(K)$ contains points reducing to~$P$.
\end{Lemma}

By scaling the variables $x$ and~$y$, we can always assume that our hyperelliptic
curve is given by a polynomial~$f$ with coefficients in the ring of integers~$\CO$
of the number field~$K$. We denote the completions of~$\CO$ and of~$K$ at the place
corresponding to a prime ideal~$\fp$ of~$\CO$ by~$\CO_\fp$ and~$K_\fp$, respectively.
We write
$F(x,z) = z^{2g+2} f(x/z)$ for the homogenization of~$f$ as an even degree binary form.
Then a smooth projective model of~$C$ is given by the weighted homogeneous equation
\[ Y^2 = F(X, Z) \]
where $X$ and~$Z$ have degree~$1$ and $Y$ has degree~$g+1$; this model sits
in a weighted projective plane.

\begin{Proposition} \label{P1}
  Let $\fp$ be a prime ideal of~$\CO$ of odd residue characteristic,
  write $k = \CO/\fp$, and set $q = \#k$. Assume
  that $q > 4g^2 - 2$. Let $f \in \CO[x]$ be squarefree (as an element of~$K[x]$) of degree $2g+2$;
  then $C \colon y^2 = f(x)$ defines a hyperelliptic curve of genus~$g \ge 1$
  over~$K$. Write $\bar{f} \in k[x]$ for the reduction mod~$\fp$ of~$f$.

  If $\bar{f}$ is not of the form $c h(x)^2$ with a polynomial $h \in k[x]$
  and a non-square $c \in k^\times$ (so in particular, $\bar{f} \neq 0$),
  then $C(K_\fp) \neq \emptyset$.
\end{Proposition}

\begin{proof}
  Let $\bar{F}$ be the reduction mod~$\fp$ of the homogenized version~$F$ of~$f$.
  By assumption, $\bar{F} \neq 0$, so we can write $\bar{F} = H^2 U$ with
  nonzero binary forms $H, U \in k[X,Z]$ and $U$ squarefree. We first consider the
  case that $U$ is not constant. Then the curve over~$k$ given by $Y^2 = U(X,Z)$
  is a nice curve of genus $\deg(U)/2 - 1$. By Lemma~\ref{L:Weil}, this curve
  has at least $q - (\deg(U) - 2) \sqrt{q} + 1$ points. At most
  $2\deg(H) = 2g + 2 - \deg(U)$ of these points satisfy $H = 0$ (there are
  at most $\deg(H)$ images $(X:Z)$ on~$\PP^1$ and each corresponds to at
  most two points), so there are at least
  \[ q - (\deg(U) - 2) \sqrt{q} + 1 - (2g + 2 - \deg(U)) > 0 \]
  \emph{smooth} $k$-points on the (usually singular) curve $Y^2 = H(X,Z)^2 U(X,Z)$.
  By Lemma~\ref{L:Hensel}, this gives us a point in~$C(K_\fp)$.
  (The inequality above follows from
  $\deg(U) - 2 \le 2 g < (q - 2g + 1)/(\sqrt{q} - 1)$.)

  If $U = c$ is constant and $c \in k^\times$ is a square, then every
  choice of $(X : Z) \in \PP^1(k)$ gives rise to $k$-points on $Y^2 = c H(X,Z)^2$,
  and there are at least $2(q + 1 - \deg(H)) = 2(q - g) > 0$ such points that are
  smooth, so that we can conclude as before. So the only remaining case
  is that $U = c$ is a constant non-square. In this case,
  $\bar{f} = c H(x,1)^2$, which is not true by assumption.
\end{proof}

So we can restrict to
\begin{enumerate}[(i)]
  \item infinite places,
  \item ``small'' odd finite places (i.e., such that $q < 4 g^2 - 2$),
  \item even finite places,
  \item odd finite places such that $\bar{f} = 0$, and
  \item odd finite places such that $\bar{f}$ is a non-square constant
        times a nonzero square.
\end{enumerate}

We discuss these sets of places in the following subsections.

A function that performs the check for local points over all completions
on hyperelliptic curves over~$\Q$ is available as~\texttt{IsLocallySolvable}
in Magma.

\subsection{Infinite places}

The infinite places are easy to deal with: we always have points at complex
places, and for a real embedding $\sigma \colon K \to \R$,
we have $C(K_\sigma) = \emptyset$ if and only if $f^\sigma$ has no real roots
and the constant term of~$f^\sigma$ is negative (here $f^\sigma \in \R[x]$ denotes
the polynomial obtained by applying~$\sigma$ to the coefficients of~$f$).

\subsection{Small odd finite places}

For the small odd finite places, we use the following procedure (which works
for any odd finite place, but can be inefficient when $q$ is large).
The case of even residue characteristic is more involved; see below.

\begin{Algorithm}[Local integral points; odd] \label{Algo:LIP}
  Let $f \in \CO[x]$ be squarefree (as an element of~$K[x]$), and let $\fp$ be a prime ideal of~$\CO$
  of odd residue characteristic with uniformizer $\pi \in \CO$.
  Let $k = \CO/\fp$ denote the residue class field;
  as before, $\bar{f}$ denotes the image of~$f$ in~$k[x]$.

  This algorithm decides whether the equation $y^2 = f(x)$ has solutions in~$\CO_\fp$.

  \begin{steplist}{1em}
    \item[0.] Set $X \leftarrow \emptyset \subset \CO[x]$.
    \item[1.] For $\xi \in k$, do the following.
              \begin{enumerate}[a.]
                \item If $\bar{f}(\xi)$ is a nonzero square, then return \textbf{true}.
                \item If $\bar{f}(\xi) = 0$ and $\bar{f}'(\xi) \neq 0$,
                      then return \textbf{true}.
                \item If $\bar{f}(\xi) = \bar{f}'(\xi) = 0$, then lift $\xi$ to $a \in \CO$; \\
                      if $f(a)$ is divisible by~$\pi^2$,
                      then set $X \leftarrow X \cup \{\pi^{-2} f(a + \pi x)\}$.
              \end{enumerate}
    \item[2.] For each $h \in X$, call this procedure recursively. \\
              If one of these calls returns \textbf{true}, then return \textbf{true}.
    \item[3.] Return \textbf{false}.
  \end{steplist}
\end{Algorithm}

To see that Algorithm~\ref{Algo:LIP} is correct, note that when the conditions
in Steps 1a or~1b are satisfied, then there is a smooth point with $x$-coordinate~$\xi$
on the reduced curve over~$k$, so by Lemma~\ref{L:Hensel}, this point
lifts to $\CO_\fp$-solutions. If in Step~1c the polynomial $f(a)$
is not divisible by~$\pi^2$, then the $\fp$-adic valuation of $f(a + \pi b)$
is equal to~$1$ for all $b \in \CO_\fp$, and so the $k$-point considered cannot lift.
Otherwise, any lift of~$(\xi,0)$ must lead to a solution of $y^2 = \pi^{-2} f(a + \pi x)$
in~$\CO_\fp$.

To see that the algorithm terminates, assume the contrary. This implies
that we have an infinite recursion. So there is a sequence~$(a_n)_{n \ge 0}$ of elements
of~$\CO$ such that for all~$n$, the polynomial
\[ f_n(x) = \pi^{-2n} f(a_0 + \pi a_1 + \pi^2 a_2 + \ldots + \pi^{n-1} a_{n-1} + \pi^n x) \]
has the property that $f_n(a_n)$ is divisible by~$\pi^2$ and $f'_n(a_n)$ is divisible by~$\pi$.
Let $a = \sum_{n=0}^\infty \pi^n a_n \in \CO_\fp$. Then, by taking limits, we
see that $f(a) = f'(a) = 0$, which contradicts the assumption that $f$ is squarefee.

\begin{Algorithm}[Small odd prime] \label{Algo:SOP}
  Let $f \in \CO[x]$ be squarefree (as an element of~$K[x]$) and of degree $2g+2$, and let $\fp$
  be a prime ideal of~$\CO$ of odd residue characteristic.
  We keep the notations $\pi$, $k$, and~$\bar{f}$ from Algorithm~\ref{Algo:LIP}.
  Recall that $F$ denotes the homogenization of~$f$ as a binary form of even degree.

  This algorithm decides whether the smooth projective curve~$C$ associated
  to $y^2 = f(x)$ has points in~$K_\fp$.

  \begin{steplist}{1em}
    \item[1.] If $\deg(\bar{f}) = 2g+2$ and the leading coefficient of~$\bar{f}$
              is a square, then return \textbf{true}.
    \item[2.] If $\deg(\bar{f}) = 2g+1$, then return \textbf{true}.
    \item[3.] If $F(1, 0)$ is divisible by~$\pi^2$ in~$\CO_\fp$, then
              call Algorithm~\ref{Algo:LIP} on~$\pi^{-2} F(1, \pi x)$. \\
              If the result is \textbf{true}, then return \textbf{true}.
    \item[4.] Call Algorithm~\ref{Algo:LIP} on~$f$. \\
              If the result is \textbf{true}, then return \textbf{true}.
    \item[5.] Return \textbf{false}.
  \end{steplist}
\end{Algorithm}

To see that Algorithm~\ref{Algo:SOP} is correct, note first that a $K_\fp$-point
on the curve~$C$ either has $x$-coordinate in~$\CO_\fp$ or the inverse
of the $x$-coordinate is in~$\pi \CO_\fp$ (this includes the case that
the point is at infinity). The existence of a point with $\fp$-adically
integral $x$-coordinate is decided in Step~4. The first two steps
check whether there is a smooth $k$-rational point at infinity; if this
is the case, then we can lift it to a $K_\fp$-rational point at infinity.
Otherwise, Step~3 checks whether $K_\fp$-rational points can reduce to
the $k$-point at infinity, and if so, call the previous algorithm to
decide whether there really are such points. The test is analogous
to Step~1c in Algorithm~\ref{Algo:LIP}.

\subsection{Even places}

The main difference compared to odd places when considering an even place
is that in characteristic~$2$, the $y$-derivative of an equation
$y^2 = f(x)$ is always zero, which is related to the fact that
a nonzero square in the residue class field does not necessarily lift
to a square in~$\CO_\fp$. This means that we need to work with a
more general form of the curve equation. The affine case is dealt with
in the following algorithm.

\begin{Algorithm}[Local integral points; even] \label{Algo:LIPE}
  Let $c \in \CO$ and $f \in \CO[x]$ be such that $4f(x) + c^2$ is squarefree
  as an element of~$K[x]$,
  and let $\fp$ be a prime ideal of~$\CO$
  of residue characteristic~$2$ with uniformizer $\pi \in \CO$.
  Let $k = \CO/\fp$ denote the residue class field;
  as before, $\bar{f}$ denotes the image of~$f$ in~$k[x]$.

  This algorithm decides whether the equation $y^2 + c y = f(x)$ has solutions in~$\CO_\fp$.

  \begin{steplist}{1em}
    \item[0.] Set $X \leftarrow \emptyset \subset \CO \times \CO[x]$.
    \item[1.] If $\bar{c} \neq 0$, then for $\xi \in k$, do the following.
              \begin{enumerate}[a.]
                \item If the equation $y^2 + \bar{c} y = \bar{f}(\xi)$ has
                      solutions in~$k$, then return \textbf{true}.
              \end{enumerate}
              Return \textbf{false}.
    \item[2.] (Now $\bar{c} = 0$.) For $\xi \in k$, do the following.
              \begin{enumerate}[a.]
                \item If $\bar{f}'(\xi) \neq 0$, then return \textbf{true}.
                \item Let $\eta \in k$ be the square root of~$\bar{f}(\xi)$.
                      Lift $\xi, \eta$ to $a, b \in \CO$. \\
                      If $\pi^2$ divides $f(a) - b^2 - c b$, then set \\
                      $X \leftarrow X \cup
                           \bigl\{\bigl(\pi^{-1}(2b+c), \pi^{-2} (f(a + \pi x) - b^2 - cb)\bigr)\bigr\}$.
              \end{enumerate}
    \item[2.] For each pair $(c', h) \in X$, call this procedure recursively. \\
              If one of these calls returns \textbf{true}, then return \textbf{true}.
    \item[3.] Return \textbf{false}.
  \end{steplist}
\end{Algorithm}

To see that Algorithm~\ref{Algo:LIPE} is correct, note that when $\bar{c} \neq 0$,
every $k$-point on the curve $y^2 + \bar{c}y = \bar{f}(x)$ is smooth
since the $y$-derivative is always $\bar{c} \neq 0$. The $x$-derivative
is $\bar{f}'(\xi)$, so the point is also smooth when that does not vanish.
In both cases, Lemma~\ref{L:Hensel} shows that there are solutions in~$\CO_\fp$.
When the $k$-points with $x$-coordinate~$\xi$ are not smooth and $f(a) - b^2 - c b$
is not divisible by~$\pi^2$, then $b^2 + c b - f(a)$ will have $\fp$-adic
valuation~$1$ for \emph{all} lifts $a, b$ to~$\CO_\fp$ of $\xi, \eta$, so no lift will
give a solution. Otherwise, we can make the indicated substition; each solution
in~$\CO_\fp$ of the original equation will give a solution to the new equation.

Termination is seen in a similar way as for Algorithm~\ref{Algo:LIP}. If the
algorithm generates an infinite recursion, then we obtain sequences $(a_n)$
and~$(b_n)$ in~$\CO$ such that with $a = \sum_{n=0}^\infty a_n \pi^n \in \CO_\fp$ and
$b = \sum_{n=0}^\infty b_n \pi^n \in \CO_\fp$, we find that $b^2 + c b = f(a)$ and
$2b + c = f'(a) = 0$, which implies that $4f(x) + c^2$ is divisible by~$(x-a)^2$,
contradicting the assumption on $f$ and~$c$.

To decide if there is a $K_\fp$-point on the projective curve, we again
separate integral and non-integral $x$-coordinates.

\begin{Algorithm}[Even prime] \label{Algo:EP}
  Let $f \in \CO[x]$ be squarefree (as an element of~$K[x]$) and of degree $2g+2$, and let $\fp$
  be a prime ideal of~$\CO$ of residue characteristic~$2$.
  We keep the notations $\pi$, $k$, and~$\bar{f}$ from Algorithm~\ref{Algo:LIPE}.
  Recall that $F$ denotes the homogenization of~$f$ as a binary form of even degree.

  This algorithm decides whether the smooth projective curve~$C$ associated
  to $y^2 = f(x)$ has points in~$K_\fp$.

  \begin{steplist}{1em}
    \item[1.] Call Algorithm~\ref{Algo:LIPE} on~$(0,f)$. \\
              If the result is \textbf{true}, then return \textbf{true}.
    \item[2.] Call Algorithm~\ref{Algo:LIPE} on~$(0, F(1, \pi x))$. \\
              If the result is \textbf{true}, then return \textbf{true}.
    \item[3.] Return \textbf{false}.
  \end{steplist}
\end{Algorithm}

\subsection{Odd places where $f$ reduces to zero}

To find the (odd) primes~$\fp$ such that $\bar{f} = 0$, we need to factor
the content of~$f$, i.e., the ideal generated by all the coefficients.
This will typically be fairly small, so the factorization step should not
be a problem.

Now assume that $\fp$ is such a prime and let $\pi \in \CO$ be a uniformizer
for~$\fp$. Let $m$ be the minimal $\fp$-adic valuation of a coefficient of~$f$;
then $m > 0$. If $m$ is even, then we can divide~$f$ by the square~$\pi^m$
to get a polynomial with nonzero reduction mod~$\fp$. Otherwise, let
$f_1 = \pi^{-m} f$; then $f$ is a square times $\pi f_1$.
If $\xi \in k$ is such that $f_1(\xi) \neq 0$, then $f(a)$ will have odd
$\fp$-adic valuation~$m$ for all $a \in \CO$ reducing to~$\xi$; so we
do not obtain points in~$C(K_\fp)$ with such $x$-coordinates. Similarly,
if $\deg(\bar{f}_1) = 2g+2$, then there are no points in~$C(K_\fp)$ with
non-integral $x$-coordinate. So we can restrict to the $\xi \in k$
such that $\bar{f}_1(\xi) = 0$ (and $\xi = \infty \in \PP^1(k)$ when
$\deg(\bar{f}_1) < 2g+2$). Fix such a~$\xi$ and lift it to $a \in \CO$.
Then we run Algorithm~\ref{Algo:LIP} on~$\pi^{-1} f_1(a + \pi x)$
(respectively, on~$\pi^{-1} F_1(1, \pi x)$ when $\xi = \infty$).

\subsection{Odd places where $f$ reduces to a non-square constant times a square}

We first need to find these primes. We can run Algorithm~\ref{Algo:SOP} for all primes
dividing both the coefficient of~$x^{2g+2}$ (which by assumption is nonzero)
and the coefficient of~$x^{2g+1}$ of~$f$. In this way, we reduce to finding
all the primes~$\fp$ that do not divide the leading coefficient~$c$ of~$f$
and have the property that $\bar{f}$ is a constant times a square (we look
at whether the constant is a square or not later).

To do this, we determine the unique monic polynomial $q \in K[x]$ of degree~$g+1$
such that $\deg(r) \le g$, where $r = f - c q^2$ and $c$ is the leading coefficient
of~$f$. The coefficients of the
polynomial~$q$ can be determined recursively from the top down to the constant
term; this involves divisions by $2$ and by~$c$, so $q, r \in \CO[(2c)^{-1}, x]$.
We factor the (numerator of the) fractional ideal of~$\CO$ generated by
the coefficients of~$r$; the prime ideals dividing it are the ones we are
looking for: if $\bar{f} = \bar{c} \tilde{q}^2$ for some monic polynomial
$\tilde{q} \in k[x]$ and $\fp \nmid 2c$, then necessarily $\tilde{q} = \bar{q}$
and therefore $\bar{r} = 0$.

We then check for each of the resulting prime ideals~$\fp$ whether $c$ reduces
to a (nonzero) square mod~$\fp$. If it does, then (recall that $\fp$ does
not divide the leading coefficient of~$f$) there are smooth points at infinity
on the reduced curve, which via Lemma~\ref{L:Hensel} implies that
$C(K_\fp) \neq \emptyset$. Otherwise, the only way
we can obtain a local point is that the $x$-coordinate reduces to a root of~$\bar{q}$.
So for each root $\xi \in k$ of~$\bar(q)$, we do the following. Let $a \in \CO$
be a lift of~$\xi$. If $\pi^2$ does not divide $f(a + \pi x)$,
then no point in~$C(K_\fp)$ can have $x$-coordinate reducing to~$\xi$.
Otherwise, we apply Algorithm~\ref{Algo:SOP} to~$\pi^{-2} f(a + \pi x)$.
If it returns~\textbf{true}, we return~\textbf{true}, otherwise we consider
the next~$\xi$. If no root~$\xi$ leads to success, we return~\textbf{false}.


\section{Rational divisors of odd degree} \label{S:Pic1}

In this section we restrict to curves of genus~$2$ over~$\Q$.
We explain how one can determine whether the curve~$C$ has
rational divisors of odd degree. This can be done whenever we know
generators of a finite index subgroup~$G_0$ of~$J(\Q)$.
In the application that is the focus of this paper, $G_0$ is generated
by the torsion subgroup~$J(\Q)_\tors$ together with the point~$P$.

Note that when $C$ has rational points, there are also rational
divisors of odd degree: every rational point gives such a divisor.
So if we can show that $C$ has no rational divisors of odd degree,
then we have shown that $C$ has no rational points.

First we saturate (see Section~\ref{S:MW}) the known subgroup~$G_0$ at~$2$.
This gives us a subgroup $G \subseteq J(\Q)$
of finite odd index; in particular $G/2G \to J(\Q)/2J(\Q)$ is an isomorphism.
Next observe that a rational divisor~$D$ of odd degree on~$C$ gives a rational point
on~$\Pic^1_C$ (by taking its class and subtracting a suitable multiple of the canonical
class~$K$, which has degree~$2$ and is defined over~$\Q$).
A rational point on~$\Pic^1_C$ does not necessarily
arise from a rational divisor, but if one of them does, then all of them do.
This is because the difference of any two rational points on~$\Pic^1_C$ is
a point in~$J(\Q)$, and all rational points on~$J$ are represented by rational
divisors (this is true more generally for hyperelliptic curves of even genus).

Then note that $\Pic^1_C$ is a $2$-covering of~$J$ (the $2$-covering map is given by
$Q \mapsto 2Q - K$), so if it has rational points, then the set of their images in~$J(\Q)$
is a coset of~$2J(\Q)$. By the isomorphism above,
we can determine a set of representatives of these cosets from~$G$.
We then test each representative if it is in the image of~$\Pic^1_C(\Q)$,
and if so, whether it can be represented by a rational divisor (and in this
case, we actually find such a divisor). If a rational divisor of odd degree
is found, we are done. Otherwise, we have shown that no rational divisor
of odd degree exists (in particular, we do not need to consider further cosets).

To check if a point in~$J(\Q)$ lifts to a rational divisor of odd degree,
we use the following diagram.
\[ \xymatrix{ \Pic^1_C \ar[r] \ar[d] & J \ar[d] \\
              \CK^\vee \ar[r]^{\delta} & \CK
            }
\]
Here $\CK$ denotes the Kummer surface $J/\{\pm 1\}$ of~$J$
and $\CK^\vee$ is the dual Kummer surface, which is the quotient
of~$\Pic^1_C$ by the involution induced by the hyperelliptic involution
on~$C$. The covering map $\Pic^1_C \to J$ descends to a twist~$\delta$
of the duplication map on the Kummer surface. Everything in the second
row of the diagram is completely explicit. So we can take a coset representative
$Q \in J(\Q)$, map it to~$\CK$ and then check if it lifts to a
rational point on~$\CK^\vee$ under~$\delta$ by solving a system
of polynomial equations (whose solution set is a zero-dimensional scheme,
so this is reasonably efficient). If a rational lift~$R$ to~$\CK^\vee$
is found, we then check if~$R$ comes from a rational divisor of odd degree.
This can be done completely analogously to the situation for genus~$3$ hyperelliptic
curves as described in~\cite{Stoll2017c}*{\S\,4}. Essentially, we check if some
expression is a nonzero square to see whether the point lifts to $\Pic^1_C(\Q)$,
and then we have to decide if the conic parameterizing the effective degree~$3$
divisors corresponding to the point on~$\Pic^1_C$ has rational points
(and find one if it does).

This functionality is available as \texttt{HasOddDegreeDivisor} in Magma.


\section{Mordell-Weil-Sieve + Chabauty} \label{S:MWS}

The basic algorithm is described in~\cite{BruinStoll2010}*{\S\,3, \S\,4.4}.
It uses the following ingredients. First of all, we fix an embedding
$i \colon C \to J$ defined over~$\Q$ (which for a curve of genus~$2$
is given in terms of a rational
point on~$C$ or an effective rational divisor of degree~$3$ on~$C$).

We begin with the ``Mordell-Weil sieve'' part. The idea originally goes back
to Scharaschkin~\cite{Scharaschkin} and was further developed by Flynn~\cite{Flynn2004};
it is as follows. We pick a finite set~$S$ of primes of good reduction for~$C$ and a
(smooth) positive integer~$N$. Consider the following commutative diagram,
where the vertical maps are given by reduction mod~$p$ for all $p \in S$.
\[ \xymatrix{ C(\Q) \ar[r]^{i} \ar@/^1.5pc/[rr]^{\phi} \ar[d] & J(\Q) \ar[r] \ar[d] &
                J(\Q)/NJ(\Q) \ar[d]^{\rho} \\
              \prod_{p \in S} C(\F_p) \ar[r]^{i} \ar@/_1.5pc/[rr]_{\psi} & \prod_{p \in S} J(\F_p) \ar[r] &
                \prod_{p \in S} J(\F_p)/NJ(\F_p)
            }
\]
We then clearly have
\begin{equation} \label{eqn:MWS}
  \phi\bigl(C(\Q)\bigr) \subseteq \rho^{-1}\bigl(\im(\psi)\bigr) \,.
\end{equation}
The set on the right is a set of cosets of~$N J(\Q)$, which can be computed explicitly
if we know that $J(\Q)_\tors + \Z \cdot P$ is saturated at all primes dividing~$N$.
See~\cite{BruinStoll2010}*{\S\,3} for a discussion how this computation can be
done efficiently (including a good choice for~$N$).
We can check the saturation condition for every prime divisor~$q$ of~$N$ and, if necessary,
replace~$P$ by some point~$P'$ such that $P = q P' + T$ with some $T \in J(\Q)_\tors$.

Now the hope is that we can use this approach to determine~$C(\Q)$. One necessary
ingrendient for this is that $C$ satisfies the following conjecture.
This is a variant of Conjecture~21 in~\cite{Stoll2011a}, which a version
of Poonen's Heuristic, restricted to cosets in~$J(\Q)$
(the original version assumes that $C(\Q)$ is empty; see~\cite{Poonen2006}).
It is also related to the Main Conjecture in~\cite{Stoll2007}.

\begin{Conjecture} \label{conj:MWS}
  For $N$ sufficiently divisible and $S$ sufficiently large depending on~$N$,
  the inclusion in~\eqref{eqn:MWS} is an equality, i.e.,
  \[ \phi\bigl(C(\Q)\bigr) = \rho^{-1}\bigl(\im(\psi)\bigr) \,. \]
\end{Conjecture}

Note that once the statement of the conjecture holds for some $N$ and~$S$, it
will stay valid for any set~$S'$ containing~$S$. (It may become invalid
when we replace $N$ by a multiple~$N'$: a coset modulo~$NJ(\Q)$ that contains
the image of a rational point will split into a number of cosets modulo~$N'J(\Q)$
that can be in~$\rho^{-1}(\im(\psi))$ but do not necessarily all contain images
of rational points.)

We also note that we can certify that the statement of Conjecture~\ref{conj:MWS}
holds for a given $N$ and~$S$: for each coset in~$\rho^{-1}(\im(\psi))$ we exhibit
a rational point on~$C$ whose image under~$\phi$ is this coset.
Conjecture~\ref{conj:MWS} ensures that we will eventually be able to do that,
by successively increasing $N$ and~$S$ until we are successful.
In this case, we know that all rational points map into a given subset of~$J(\Q)/NJ(\Q)$,
but it is still possible that some of the fibers of~$\phi$ contain more than one point.

Since we know that $C(\Q)$ is finite and $C(\Q)$ injects (via~$i$) into~$J(\Q)$,
which is a finitely generated abelian group, it follows that for $N$ sufficiently
divisible, $\phi$ must be injective. If this holds for~$N$, then our algorithm
will have determined~$C(\Q)$. The problem we have in practice is that we need
an explicit criterion that allows us to verify that $\phi$ is injective for a given~$N$.

This is where the ``Chabauty'' part comes in. It is based on the following fact.
Let $p$ be a prime, which we assume to be of good reduction for~$C$ for simplicity.
There is a pairing (see~\cite{Coleman,McCP,Stoll2006a})
\[ J(\Q_p) \times \Omega^1_{J/\Q_p} \To \Q_p\,, \qquad
   (P, \omega) \longmapsto \int_0^P \omega = \omega(\log P) \,,
\]
where $\log \colon J(\Q_p) \to T_0 J(\Q_p)$ is the $p$-adic abelian logarithm
and $\omega$ is identified with an element of the cotangent space~$T^*_0 J(\Q_p)$.
Since $J(\Q)$ has rank~$1$ by assumption and $\Omega^1_{J/\Q_p}$ is a $2$-dimensional
vector space over~$\Q_p$, there must be a non-zero differential \hbox{$\omega_p \in \Omega^1_{J/\Q_p}$}
that annihilates~$J(\Q)$ under this pairing. We call such an~$\omega_p$ an
\emph{annihilating differential}. There is a canonical isomorphism between the
spaces of regular $1$-forms on~$J$ and on~$C$, so we can consider an
annihilating differential~$\omega_p$ as a regular differential on~$C_{\Q_p}$.
It then follows that for all pairs $Q, Q' \in C(\Q)$ of rational points,
\[ \int_Q^{Q'} \omega_p = 0 \,. \]
When $Q$ and~$Q'$ reduce to the same point in~$C(\F_p)$, this integral can be
computed by evaluating the formal integral of a power series representing~$\omega_p$.
We can scale~$\omega_p$ such that its reduction~$\bar{\omega}_p$ mod~$p$ makes sense
and is nonzero.
By~\cite{Stoll2006a}*{Prop.~6.3}, the number of rational points on~$C$ that
reduce to $\bar{Q} \in C(\F_p)$ is at most $1$ plus the order of vanishing
of~$\bar{\omega}_p$ at~$\bar{Q}$ (unless $p$ is very small). In particular,
there can be at most one rational point reducing to~$\bar{Q}$ when $\bar{\omega}_p$
does not vanish at~$\bar{Q}$ and $p > 2$ (by~\cite{Stoll2006a}*{Lemma~6.1},
we have $\delta(p, 0) = 0$ for $p > 2$). This shows that when $p > 2$
and $\bar{\omega}_p$ does not vanish at any point in~$C(\F_p)$, then
the reduction map $C(\Q) \to C(\F_p)$ is injective. This implies that the map~$\phi$
is injective when $N$ is a multiple of the exponent of (the image of~$J(\Q)$ in)
$J(\F_p)$. (Note that $\bar{\omega}_p$ vanishes at the points with a certain
$x$-coordinate $\xi \in \F_p$; the condition is that their $y$-coordinates are
not in~$\F_p$, i.e., $\bar{f}(\xi)$ is a non-square in~$\F_p$.)

So we want to pick~$N$ in such a way that it is a multiple of the exponent
of~$J(\F_p)$ for a good prime $p > 2$ such that $\bar{\omega}_p$ does not vanish
on~$C(\F_p)$. This means that we have to find at least one such prime.
The following conjecture says that we will easily find such primes unless
there is a good reason why this is impossible. Such a good reason is provided
by the geometry of~$J$: if $J$ \emph{splits}, which means that it is isogenous
over~$\Q$ to the product of two elliptic curves $E$ and~$E'$, then one of
the elliptic curves, say~$E$, must have rank~$1$ and the other one must have
rank~$0$. The pull-back~$\omega$ of a nonzero regular differential on~$E'$ then vanishes
along the image of~$E$ in~$J$ (obtained via Picard functoriality from
the dominant morphism $C \stackrel{i}{\to} J \to E \times E' \to E$) and
will therefore be an annihilating differential for all (good) primes~$p$.

\begin{Conjecture} \label{conj:diff}
  The set of primes~$p$ of good reduction for~$C$ such that the reduction
  mod~$p$, $\bar{\omega}_p$, of a suitably scaled annihilating differential~$\omega_p$
  does not vanish on~$C(\F_p)$ has positive density, unless $J$ splits and
  the associated global annihilating differential~$\omega$ vanishes at a
  rational point of~$C$.
\end{Conjecture}

This is a strengthening of Conjecture~4.2 in~\cite{BruinStoll2010} in the
case when $C$ has genus~$2$, taking into account the heuristic computation
of the density, which for hyperelliptic curves (and ``rank defect'' $g - r = 1$)
predicts density~$1/2$. The conjecture in loc.~cit.\ makes a statement
for simple, so non-split, Jacobians. If in the split case, $\omega$ does not
vanish at a rational point, then it vanishes at a pair of points with
$x$-coordinate $\xi \in \Q$ such that $f(\xi)$ is a non-square in~$\Q$
(or $\xi = \infty$ and the leading coefficient of~$f$ is a non-square).
There is then a set of primes of density~$1/2$ such that $f(\xi)$ is a
non-square mod~$p$; for such~$p$, $\bar{\omega}_p = \bar{\omega}$ will
not vanish on~$C(\F_p)$.

So, assuming the conjecture, we are in good shape \emph{unless} $J$ splits
and $\omega$ vanishes at a rational point. In this case, the curve~$C$
has morphisms of some degree~$d$ (which we can assume to be minimal) to
the elliptic curves~$E$ and~$E'$, with $E(\Q)$ of rank~$1$ and~$E'(\Q)$ finite.
If we can determine the morphism $C \to E'$ explicitly, the we can solve
our problem by first determining the finitely many points in~$E'(\Q)$ and
then determining the rational points on~$C$ in each fiber above one of
these points.

There are fairly simple concrete criteria for when such a splitting of~$J$ exists
with $d = 2$ or $d = 3$, and in this case there are formulas for the elliptic
curves and the morphisms. This is implemented via the functions
\texttt{Degree2Subcovers} and~\texttt{Degree3Subcovers} in Magma.
However, larger degrees~$d$ do occur, as the following examples show.

\begin{Examples}
  Here are some examples of genus~$2$ curves~$C$ that have maps to elliptic
  curves of minimal degree $d > 3$.
  \begin{enumerate}[(1)] \setcounter{enumi}{3}
    \item $y^2 = 6 x^5 + 28 x^3 + 54 x$ with $d = 4$.
    \item $y^2 = 64 x^6 + 180 x^3 + 125$ with $d = 5$.
    \item $y^2 = 192 x^5 + 420 x^4 + 504 x^3 + 177 x^2 + 66 x + 9$ with $d = 6$.
    \item $y^2 = 4 x^6 - 12 x^5 + 81 x^4 - 22 x^3 + 181 x^2 + 808 x + 304$ with $d = 7$.
    \item $y^2 = x^6 - 6 x^5 + 23 x^4 - 32 x^3 + 71 x^2 + 126 x + 213$ with $d = 8$.
  \end{enumerate}
  These were found among a set of about $6$~million curves in a list compiled
  by Drew Sutherland; they caused the previous version of the \texttt{Chabauty}
  procedure in Magma to get stuck.
\end{Examples}

\begin{Remark}
  The two elliptic curves in the product isogenous to a split genus~$2$ Jacobian
  have isomorphic $d$-torsion Galois modules; such elliptic curves are said
  to be \emph{$d$-congruent}.
  Conversely, from a pair of $d$-congruent elliptic curves (such that the
  isomorphism of $d$-torsion modules is not induced by an isomorphism
  of the curves), one obtains a split genus~$2$ Jacobian.
  There are split genus~$2$ Jacobians for which the degree~$d$ is even larger
  than in the examples above.
  See work of Fisher~\cite{Fisher2014,Fisher2015,Fisher2018,Fisher2019,Fisher2021}
  on congruent elliptic curves.
  Fisher~\cite{Fisher2021} gives two examples with $d = 17$
  and conjectures (Conjecture~1.1 in \emph{loc.~cit.}) that for primes $d \ge 17$,
  apart from examples coming from isogenies,
  these are the only examples up to quadratic twist, whereas~\cite{Fisher2019}
  gives an infinite family of pairs of $13$-congruent elliptic curves.
  Fisher's conjecture is a strong form of the Frey-Mazur Conjecture, which states
  that for all sufficiently large primes~$p$, all $p$-congruences of elliptic curves
  come from isogenies. (Note that an isogeny $E \to E'$ of degree prime to~$p$
  induces an isomorphism $E[p] \cong E'[p]$ of Galois modules.)
\end{Remark}

Luckily, it turns out that it is not actually necessary to find the morphism $C \to E'$
explicitly. Recall that we always assume that the rank of~$J(\Q)$ is~$1$.

\begin{Lemma} \label{L:split}
  Assume that there is a non-constant morphism $C \to E$ with an elliptic curve~$E$.
  Let $Q \in C(\Q)$ be such that the global annihilating differential~$\omega \in \Omega^1_C(\Q)$
  vanishes at~$Q$. Let $p > 3$ be a prime of good reduction for~$C$.
  Assume that either $Q$ is a Weierstrass point or that the reduction $\bar{Q} \in C(\F_p)$
  is not a Weierstrass point. Then the map $C(\Q) \to C(\F_p)$ is injective.
\end{Lemma}

\begin{proof}
  We already know that the fibers of $C(\Q) \to C(\F_p)$ have at most one element
  above points in~$C(\F_p)$ at which $\bar{\omega}$ does not vanish.
  It remains to show that the fiber above the reduction~$\bar{Q}$ of~$Q$ mod~$p$
  consists of the single point~$Q$. (The same argument will work for the other
  rational point~$\iota(Q)$ on which $\omega$ vanishes, when $Q$ is not a Weierstrass point.)

  By the definition of an annihilating differential, for any point $Q' \in C(\Q)$
  such that $\bar{Q}' = \bar{Q}$ the $p$-adic integral
  \[ \int_Q^{Q'} \omega \]
  vanishes. This integral can be computed by formally integrating a power series
  and evaluating at a suitable parameter (it is a ``tiny integral'').

  First assume that $Q$ is not a Weierstrass point and its reduction~$\bar{Q} \in C(\F_p)$
  is also not a Weierstrass point (i.e., the numerator of the $y$-coordinate of~$Q$ is
  not divisible by~$p$). Then $t = x - x(Q)$ when $x(Q) \in \Z_p$ or
  $t = x^{-1} - x(Q)^{-1}$ when $x(Q)^{-1} \in \Z_p$ is a local parameter at~$Q$
  reducing to a local parameter at~$\bar{Q}$. Since $\omega$ vanishes at~$Q$, we have
  (after possibly scaling~$\omega$)
  \[ \omega = t w(t)\,dt \]
  with a power series $w(t) = 1 + a_1 t + a_2 t^ 2 + \ldots$ with coefficients in~$\Z_p$,
  and the integral from $Q$ to $Q'$ (in the same $p$-adic residue disk) is
  \[ \int_Q^{Q'} \omega = \int_0^{t(Q')} (t + a_1 t^2 + ...) \, dt
                       = \frac{1}{2} t(Q')^2 + \frac{a_1}{3} t(Q')^3 + \frac{a_2}{4} t(Q')^4 + \ldots \,;
  \]
  this has a double zero at~$Q$ and (since $p > 3$) no further zeros on the residue disk.

  If $Q$ is a Weierstrass point, then $t = y$ when $Q \neq \infty$ or $t = y/x^3$ when $Q = \infty$
  is again a local parameter at~$Q$ reducing to a local parameter at~$\bar{Q}$.
  Since $\omega$ in this case vanishes to order~$2$ at~$Q$, we have (again after scaling)
  \[ \omega = t^2 w(t)\,dt \]
  with an even power series $w(t) = 1 + a_2 t^2 + a_4 t^4 + \ldots$, and
  \[ \int_Q^{Q'} \omega = \int_0^{t(Q')} (t^2 + a_2 t^4 + \ldots) \, dt
                        = \frac{1}{3} t(Q')^3 + \frac{a_2}{5} t(Q')^5 + \ldots \,;
  \]
  this has a triple zero at~$Q$ and (again since $p > 3$, note the increase by~$2$
  of the degree) no further zeros on the residue disk.
\end{proof}

Note that when $\omega$ vanishes at a non-Weierstrass point~$Q$ that reduces mod~$p$
to a Weierstrass point~$\bar{Q}$, then the fiber of the reduction map above~$\bar{Q}$
contains (at least) the two rational points $Q$ and $\iota(Q)$ (where $\iota$ denotes
the hyperelliptic involution), so the conclusion cannot hold in this case.
But in any case, we see that the reduction map will be injective for all but finitely
many primes, with an explicit set of possible exceptions.

So if we seem unable to find a good odd prime~$p$ such that $\bar{\omega}_p$ does not
vanish on~$C(\F_p)$, we try to first show that $J$ splits (and if so, determine the
degree~$d$). This can be done numerically by computing a period matrix
(Magma contains functionality for analytic Jacobians of hyperelliptic curves;
see~\cite{CMSV}).
The result we obtain is not rigorous, but it will be certified by the next step.
This uses the following observation.

\begin{Lemma}
  Let $C \to E$ be a morphism of degree~$d$ and let $\tilde{E} \subset J$
  be the image of~$E$ under the induced morphism $E \to J$. Let $\CK$ be the
  Kummer surface of~$J$ and denote by~$Y$ the image of~$\tilde{E}$ on~$\CK$.
  Then $\tilde{E} \to Y$ is a double cover ramified in four points, so
  $Y$ is a smooth rational curve of degree~$d$ in~$\PP^3$. Also, $Y$
  is contained in a surface of degree~$\lceil d/2 \rceil$ that does not
  contain~$\CK$.
\end{Lemma}

\begin{proof}
  Since $E \to J$ is a homomorphism of abelian varieties and $J \to \CK$
  identifies points with their negatives, the morphism $E \to J$ descends
  to a morphism $\PP^1 \to \CK$ (where $E \to \PP^1$ is the $x$-coordinate map);
  in particular, $Y$ is a rational curve. The degree of~$\tilde{E}$ with respect
  to the theta divisor on~$J$ is~$d$; this implies that the degree of~$Y$
  is also~$d$. We have that $J[2] \cap \tilde{E} = \tilde{E}[2]$, which is of
  size~$4$; as $J[2]$ is the set of ramification points of $J \to \CK$, this
  shows that $\tilde{E} \to Y$ is ramified in four points (this also follows
  from Riemann-Hurwitz).

  The morphism $\PP^1 \to Y \subset \CK \subset \PP^3$ is given by a quadruple
  of binary forms of degree~$d$. There is a $\binom{n+3}{3}$-dimensional
  space of degree~$n$ forms in four variables containing an $\binom{n-1}{3}$-dimensional
  subspace of multiples of the defining quartic equation of~$\CK$.
  So there will be a degree~$n$ form vanishing on~$Y$ that is not a multiple
  of the defining equation of~$\CK$ whenever
  \[ 2n^2 + 2 = \binom{n+3}{3} - \binom{n-1}{3} > dn + 1 \,, \]
  where $dn + 1$ is the dimension of the space of binary forms of degree~$dn$.
  This is the case for $n \ge \lceil d/2 \rceil$.
\end{proof}

The idea now is to find~$Y$ by interpolating points. Note that $\tilde{E}$
will contain a point of the form $P' = n P + T$ for some $n \ge 1$ and $T \in n J(\Q)_\tors$,
where $n$ divides the least common multiple of~$d$ and the exponent of~$J(\Q)_\tors$.
Since $\tilde{E}$ is an abelian subvariety, it will then contain all multiples of~$P'$,
so $Y$ will contain all their images on~$\CK$.

For $n = 1$, $n = d$, $n = \operatorname{lcm}(d, \exp(J(\Q)_\tors))$
and in each case for all possibilities of $T \in n J(\Q)_\tors$,
we therefore compute enough of the multiples of $P' = nP + T$ and find all hypersurfaces
of degree $m = \lceil d/2 \rceil$ that pass through these points.
If there are any such hypersurfaces, we then check for whether
their intersection with the Kummer surface
has dimension~$1$ and contains a curve of degree~$d$ and genus~$0$ passing through
exactly four points of order~$2$ and containing the image of~$P'$.
If this is the case for some choice of $n$ and~$T$ as above, then we have shown
that $J$ splits and we have found
the image~$Y$ of~$\tilde{E}$ on the Kummer surface. We can then obtain the $x$-coordinate
of the points on which the annihilating differential vanishes from the
tangent line of~$Y$ at the origin.
It remains to verify that there are rational points on~$C$ with that $x$-coordinate.
If this is the case, then we can apply Lemma~\ref{L:split} above.

Magma's \texttt{Chabauty} function carries out the Mordell-Weil Sieve +
Chabauty computation explained in this section. It includes the check for
rational divisors of odd degree. More precisely, it first searches for
rational points on the curve (unless a rational point on~$C$ is given as an
optional argument; then that is used to embed $C$ into~$J$), and if it does
not find one up to some height bound, then it runs \texttt{HasOddDegreeDivisor}.
The search for rational points is done using \texttt{ratpoints}~\cite{StollRP};
this code is part of~Magma.


\section{Concluding remarks.}

The functionality described here (except the part that is summarized in Section~\ref{S:MW},
which was contributed in~2019/2020) forms part of what the new (as of version~2.29)
Magma function \texttt{RationalPointsGenus2} is doing. Its purpose is to try
to provably determine the full set of rational points on a curve~$C$ of genus~$2$
over~$\Q$. It essentially performs the following steps.

\begin{enumerate}[1.]
  \item Search for small rational points.
  \item If no point was found, then test whether $C$ has points everywhere locally
        (see Section~\ref{S:ELS}). If $C$ fails to have points over some
        completion of~$\Q$, then $C(\Q) = \emptyset$; Stop.
  \item Search again for rational points up to a larger bound.
  \item If no point was found, then compute the (fake) $2$-Selmer set of~$C$;
        see~\cite{BruinStoll2009}. If this is empty, then $C(\Q) = \emptyset$; Stop.
  \item Check if $C$ admits morphisms of degree $2$ or~$3$ to elliptic curves.
        If so, and one of the elliptic curves, say~$E$, has finite group of rational points,
        then determine $C(\Q)$ by finding the rational points in the fibers above
        the rational points of~$E$; Stop.
  \item Attempt to determine the rank~$r$ of~$J(\Q)$ (see Section~\ref{S:MW}).
        If not successful, then report \emph{Failure;} Stop.
  \item If $r = 0$, then determine $C(\Q)$ by pulling back the finitely many torsion
        points in~$J(\Q)$ under $C \to J$, $P \mapsto [P - \iota(P)]$; Stop.
  \item Now $r \ge 1$. If no rational point was found so far, then test
        whether $C$ has rational divisors of odd degree (see Section~\ref{S:Pic1}).
        If this is not the case, then $C(\Q) = \emptyset$; Stop.
  \item If $r = 1$, then run the Mordell-Weil Sieve + Chabauty computation
        (see Section~\ref{S:MWS}). If this terminates, it will
        have determined~$C(\Q)$; Stop. (In the actual implementation the computation
        will be aborted and \emph{Failure} is reported when certain bounds are reached.)
  \item Now $r \ge 2$. If no rational point was found so far, use the rational
        divisor of odd degree that was found earlier to obtain an embedding $C \to J$ and
        run a Mordell-Weil Sieve computation as in~\cite{BruinStoll2010}.
        If this proves that $C(\Q) = \emptyset$,
        then return $\emptyset$, otherwise report \emph{Failure;} Stop.
  \item Report \emph{Failure;} Stop.
\end{enumerate}

Assuming Conjectures \ref{conj:MWS} and~\ref{conj:diff}, this algorithm will succed
in determining~$C(\Q)$ whenever one of the following conditions holds.
\begin{enumerate}[(1)]
  \item The $2$-Selmer set of~$C$ is empty (this includes the case that
        $C$ fails to have points over all completions of~$\Q$).
  \item There is an elliptic subcover~$E$ of degree $2$ or~$3$, for which
        it can be determined that the rank is zero.
  \item The rank~$r$ of~$J(\Q)$ can be determined and one of the following holds.
        \begin{enumerate}[(a)]
           \item There is no rational divisor of odd degree on~$C$.
           \item $r \le 1$.
           \item $C(\Q) = \emptyset$.
        \end{enumerate}
\end{enumerate}

Note that Conjecture~\ref{conj:MWS} implies that the Mordell-Weil Sieve computation
will detect that $C(\Q) = \emptyset$ whenever this is the case.

Drew Sutherland has constructed a database of more than 6~million curves of
genus~$2$ over~$\Q$ of conductor up to~$2^{20}$ \cite{Sutherland},
which extends the database~\cite{BSSVY}
currently available in the LMFDB~\cite{lmfdb}. \texttt{RationalPointsGenus2}
successfully determines the set of rational points on about~92\% of these
curves (assuming GRH for the rank computation).
This high success rate is probably explained to a large part by the fact
that most curves do not have rational points: heuristically, one would
expect only a fraction of~$\ll N^{-1/2}$ curves $y^2 = f(x)$ with integral
coefficients bounded by~$N$ in absolute value to have rational points.
(In the data set, $2\,528\,131$ of $6\,216\,959$~curves have a known rational point.)
Experimentally, a fairly large proportion of genus~$2$ curves without rational points
have an empty $2$-Selmer set (see~\cite{BruinStoll2009}*{Section~10}),
so in these cases the algorithm will be successful
even without attempting to determine the rank of~$J(\Q)$. Among the curves
for which the rank can be determined, we expect a large majority to have
rank $0$ or~$1$, in which cases the algorithm will also be successful.

We remark that there is currently no generally applicable (and practical) algorithm
available that can determine~$C(\Q)$ when $r \ge 2$ and $C(\Q) \neq \emptyset$.
However, there are some restricted cases when this is possible. One such case
is when $J$ has Néron-Severi rank~$\rho$ strictly larger than~$1$; this occurs when
$J$ has real multiplication or is split. When (for a general curve of genus~$g$)
$r < g - 1 + \rho$, then \emph{Quadratic Chabauty}
\cite{BalakrishnanDogra2018,BDMTV,BalakrishnanDogra2021,BBBM,EdixhovenLido}
can determine~$C(\Q)$. There are also recent steps~\cite{Dogra2023,Dogra2024}
in the direction of making
the full second level of the Chabauty-Kim approach effective and feasible;
this would lead to an algorithm that can determine~$C(\Q)$ when $r \le g^2$.

Another observation is that every rational point on~$C$ lifts to one of finitely
many $2$-covering curves~$D_\xi$ of genus~$17$, where $\xi$ runs through the
elements of the $2$-Selmer set of~$C$. These curves~$D_\xi$ each have $15$
maps to elliptic curves that are in general not defined over~$\Q$, however;
the Galois action on them corresponds to the Galois action on the points
of order~$2$ on~$J$, which in turn correspond to the factorizations of~$f$
into a quadric and a quartic (up to scaling). If, for each~$\xi$, there is
one such map $D_\xi \to E_\xi$, which is defined over a number field~$K_\xi$
of degree~$d$ and such that the rank of~$E_\xi(K_\xi)$ is strictly less than~$d$,
then a variant of Chabauty's method, known as \emph{Elliptic Curve Chabauty}
can be used; see~\cite{Bruin2003}. In practice, this requires the degrees~$d$
to be reasonably small (say, $d \le 5$ or so); otherwise, the determination
of~$E_\xi(K_\xi)$ is likely infeasible. This means that the Galois group of~$f$
has to be quite small. So this approach is practical only in fairly limited
situations. In any case, implementing a version of this and including it
in \texttt{RationalPointsGenus2} is a project for the not-too-distant future.


\section*{Acknowledgments}

This work was supported by a grant from the Simons Foundation International
[SFI-MPS-Infrastructure-00008651, AVS] that provided funding for a visit
of the author to MIT in March~2025, during which a major part of the work described
in this paper was done.


\end{document}